\input amstex
\documentstyle{amsppt}
%----------------------------------------------------------------
% Title:     A note on rational and elliptic curves associated 
%            with the cuboid factor equations.
% Author:    Ruslan Sharipov
% Comments:  AmSTeX, 15 pages, amsppt style
% MSC-class: 11D25, 11D72, 14E05, 14H45, 14H52 
%----------------------------------------------------------------
%           Replacement for output macro definition
%
\catcode`@=11
\redefine\output@{%
  \def\break{\penalty-\@M}\let\par\endgraf
  \ifodd\pageno\global\hoffset=105pt\else\global\hoffset=8pt\fi  
  \shipout\vbox{%
    \ifplain@
      \let\makeheadline\relax \let\makefootline\relax
    \else
      \iffirstpage@ \global\firstpage@false
        \let\rightheadline\frheadline
        \let\leftheadline\flheadline
      \else
        \ifrunheads@ %\let\makefootline\relax
        \else \let\makeheadline\relax
        \fi
      \fi
    \fi
    \makeheadline \pagebody \makefootline}%
  \advancepageno \ifnum\outputpenalty>-\@MM\else\dosupereject\fi
}
\def\Beta{\mathchar"0\hexnumber@\rmfam 42}
\catcode`\@=\active
%----------------------------------------------------------------
\nopagenumbers
\chardef\textvolna='176
\font\twelvesmmi=cmmi12
\chardef\bigalpha='013
\def\negskp{\hskip -2pt}

\chardef\degree="5E
%\def\compos{\,\raise 1pt\hbox{$\sssize\circ$} \,}
%\def\id{\operatorname{id}}
%\font\eightrm=cmr8
%\def\LT{\operatorname{\text{\eightrm LT}}}
%\def\LM{\operatorname{\text{\eightrm LM}}}
%\def\LC{\operatorname{\text{\eightrm LC}}}
%\accentedsymbol\hatgamma{\kern 2pt\hat{\kern -2pt\gamma}}
%\accentedsymbol\checkgamma{\kern 2.5pt\check{\kern -2.5pt\gamma}}
\def\blue#1{#1}

\catcode`#=11\def\diez{#}\catcode`#=6
\catcode`&=11\catcode`&=4
\catcode`_=11\def\podcherkivanie{_}\catcode`_=8
\catcode`~=11\def\volna{~}\catcode`~=\active
\def\mycite#1{\cite{\blue{#1}}\immediate\special{ps:
     ShrHPSdict begin /ShrBORDERthickness 0 def}}
\def\myciterange#1#2#3#4{\cite{\blue{#2#3#4}}\immediate\special{ps:
     ShrHPSdict begin /ShrBORDERthickness 0 def}}
\def\mytag#1{%
    \tag#1}
\def\mythetag#1{\thetag{\blue{#1}}\immediate\special{ps:
     ShrHPSdict begin /ShrBORDERthickness 0 def}}
\def\myrefno#1{\no#1}
\def\myhref#1#2{\blue{#2}\immediate\special{ps:
     ShrHPSdict begin /ShrBORDERthickness 0 def}}
\def\myEarXivlink{\myhref{http://arXiv.org}{http:/\negskp/arXiv.org}}

\def\mytheorem#1{\csname proclaim\endcsname{Theorem #1}}
\def\mytheoremwithtitle#1#2{\csname proclaim\endcsname{Theorem #1#2}}
\def\mythetheorem#1{\blue{#1}\immediate\special{ps:
     ShrHPSdict begin /ShrBORDERthickness 0 def}}
\def\mylemma#1{\csname proclaim\endcsname{Lemma #1}}
\def\mylemmawithtitle#1#2{\csname proclaim\endcsname{Lemma #1#2}}
\def\mythelemma#1{\blue{#1}\immediate\special{ps:
     ShrHPSdict begin /ShrBORDERthickness 0 def}}
\def\mycorollary#1{\csname proclaim\endcsname{Corollary #1}}

\def\mydefinition#1{\definition{Definition #1}}

\def\myconjecture#1{\csname proclaim\endcsname{Conjecture #1}}
\def\myconjecturewithtitle#1#2{\csname proclaim\endcsname{Conjecture #1#2}}

\def\myproblem#1{\csname proclaim\endcsname{Problem #1}}
\def\myproblemwithtitle#1#2{\csname proclaim\endcsname{Problem #1#2}}
\def\mytheproblem#1{\blue{#1}\immediate\special{ps:
     ShrHPSdict begin /ShrBORDERthickness 0 def}}

%----------------------------------------------------------------
% Cyrillic fonts definition
\font\eightcyr=wncyr8
%\font\tencyr=wncyr10
%----------------------------------------------------------------
\pagewidth{360pt}
\pageheight{606pt}
\topmatter
\title
A note on rational and elliptic curves associated 
with the cuboid factor equations.
\endtitle
\rightheadtext{A note on rational and elliptic curves \dots}
\author
Ruslan Sharipov
\endauthor
\address Bashkir State University, 32 Zaki Validi street, 450074 Ufa, Russia
\endaddress
\email\myhref{mailto:r-sharipov\@mail.ru}{r-sharipov\@mail.ru}
\endemail
\abstract
     A rational perfect cuboid is a rectangular parallelepiped whose edges and 
face diagonals are given by rational numbers and whose space diagonal is equal 
to unity. It is described by a system of four equations with respect to six 
variables. The cuboid factor equations were derived from these four equations 
by symmetrization procedure. They constitute a system of eight polynomial 
equations, which has been solved parametrically. In the present paper its
parametric solution is expressed through intersections or rational and elliptic 
curves arranged into parametric families. 
\endabstract
\endtopmatter
%\loadbold
%\loadeufb
\TagsOnRight
\document
% \input countstyle

%\special{header=resource.eps}
\head
1. Introduction.
\endhead
     Finding a rational perfect cuboid is equivalent to finding a perfect
cuboid with all integer edges and diagonals, which is an old unsolved problem
(see \myciterange{1}{1}{--}{44}). Here are the equations describing perfect 
cuboids:
$$
\xalignat 2
&\hskip -2em
p_{\kern 1pt 0}=x_1^2+x_2^2+x_3^2-L^2=0,
&&p_{\kern 1pt 1}=x_2^2+x_3^2-d_1^{\kern 1pt 2}=0,\\
\vspace{-1.7ex}
\mytag{1.1}\\
\vspace{-1.7ex}
&\hskip -2em
p_{\kern 1pt 2}=x_3^2+x_1^2-d_2^{\kern 1pt 2}=0,
&&p_{\kern 1pt 3}=x_1^2+x_2^2-d_3^{\kern 1pt 2}=0.
\endxalignat
$$
The variables $x_1$, $x_2$, $x_3$ in \mythetag{1.1} represent three edges of a 
cuboid, the variables $d_1$, $d_2$, $d_3$ are its face diagonals, and $L$ is its 
space diagonal. In the case of a rational perfect cuboid we set $L=1$. Using the 
polynomials $p_{\kern 1pt 0}$, $p_{\kern 1pt 1}$, $p_{\kern 1pt 2}$, and 
$p_{\kern 1pt 3}$ from \mythetag{1.1}, we write the following eight polynomial
equations:
$$
\xalignat 2
&\hskip -2em
p_{\kern 1pt 0}=0,
&&\sum^3_{i=1}p_{\kern 1pt i}=0,\\
\vspace{1ex}
&\hskip -2em
\sum^3_{i=1}d_i\,p_{\kern 1pt i}=0,
&&\sum^3_{i=1}x_i\,p_{\kern 1pt i}=0,\\
\vspace{-1.7ex}
\mytag{1.2}\\
\vspace{-1.7ex}
&\hskip -2em
\sum^3_{i=1}d_i^{\kern 1pt 2}\,p_{\kern 1pt i}=0,
&&\sum^3_{i=1}x_i^2\,p_{\kern 1pt i}=0,\\
\vspace{1ex}
&\hskip -2em
\sum^3_{i=1}x_i\,d_i\,p_{\kern 1pt i}=0,
&&\sum^3_{i=1}x_i^2\,d_i^{\kern 1pt 2}\,p_{\kern 1pt i}=0.\\
\endxalignat
$$
The equations \mythetag{1.2} are called the cuboid factor equations. They were 
derived from the original cuboid equations \mythetag{1.1} as a result of a 
symmetry approach initiated in \mycite{45} (see also \myciterange{46}{46}{--}{48}).
\par
     It is easy to see that each solution of the original cuboid equations 
\mythetag{1.1} is a solution for the factor equations \mythetag{1.2}. Generally 
speaking, the converse is not true. However, in \mycite{47} the following theorem 
was proved. 
\mytheorem{1.1} Each integer or rational solution of the cuboid factor equations 
\mythetag{1.2} such that $x_1>0$, $x_2>0$, $x_3>0$, $d_1>0$, $d_2>0$, and $d_3>0$ 
is an integer or rational solution for the original cuboid equations \mythetag{1.1}.
\endproclaim
      Due to Theorem~\mythetheorem{1.1} the cuboid factor equations \mythetag{1.2} 
are equivalent to the equations \mythetag{1.1} in studying perfect cuboids. But in 
this paper, saying a solution of the factor equations, we assume any integer or 
rational solution, i\.\,e\. even such that some of the inequalities $x_1>0$, $x_2>0$, 
$x_3>0$, $d_1>0$, $d_2>0$, $d_3>0$ or all of them are not fulfilled.\par 
      The left hand sides of the cuboid factor equations are multisymmetric
polynomials in $x_1$, $x_2$, $x_3$ and $d_1$, $d_2$, $d_3$, i\.\,e\. they are 
invariant with respect to the permutation group $S_3$ acting upon $x_1$, $x_2$, $x_3$, 
$d_1$, $d_2$, $d_3$, and $L$ as follows:
$$
\xalignat 3
&\sigma(x_i)=x_{\sigma i},
&&\sigma(d_i)=d_{\sigma i},
&&\sigma(L)=L.
\endxalignat
$$ 
For the theory of multisymmetric polynomials the reader is referred to
\myciterange{49}{49}{--}{69}. According to this theory, each multisymmetric
polynomial of $x_1$, $x_2$, $x_3$ and $d_1$, $d_2$, $d_3$ is expressed through 
the following nine elementary multisymmetric polynomials:
$$
\gather
\hskip -2em
\aligned
&x_1+x_2+x_3=E_{10},\\
&x_1\,x_2+x_2\,x_3+x_3\,x_1=E_{20},\\
&x_1\,x_2\,x_3=E_{30},
\endaligned
\mytag{1.3}\\
\vspace{2ex}
\hskip -2em
\aligned
&d_1+d_2+d_3=E_{01},\\
&d_1\,d_2+d_2\,d_3+d_3\,d_1=E_{02},\\
&d_1\,d_2\,d_3=E_{03},
\endaligned
\quad
\mytag{1.4}\\
\vspace{2ex}
\hskip -2em
\aligned
&x_1\,x_2\,d_3+x_2\,x_3\,d_1+x_3\,x_1\,d_2=E_{21},\\
&x_1\,d_2+d_1\,x_2+x_2\,d_3+d_2\,x_3+x_3\,d_1+d_3\,x_1=E_{11},\\
&x_1\,d_2\,d_3+x_2\,d_3\,d_1+x_3\,d_1\,d_2=E_{12}.
\endaligned
\mytag{1.5}
\endgather
$$
Expressing the left hand sides of the cuboid factor equations \mythetag{1.2} 
through the polynomials \mythetag{1.3}, \mythetag{1.4}, and \mythetag{1.5}, 
one gets polynomial equations with respect to the variables $E_{10}$, $E_{20}$, 
$E_{30}$, $E_{01}$, $E_{02}$, $E_{03}$, $E_{21}$, $E_{11}$, $E_{12}$, and $L$ 
(see \thetag{3.1} through \thetag{3.7} in \mycite{48}). These equations were 
complemented with fourteen identities expressing the mutual algebraic dependence 
of the elementary multisymmetric polynomials \mythetag{1.3}, \mythetag{1.4},
and \mythetag{1.5} (see \thetag{3.8} in \mycite{48}). As a result a system of 
twenty two equations was obtained. In \mycite{48} this huge system of 
twenty two equations was reduced to the following single polynomial equation for 
$E_{10}$, $E_{01}$, $E_{11}$, and $L$:
$$
\hskip -2em
(2\,E_{11})^2+(E_{01}^2+L^2-E_{10}^2)^2-8\,E_{01}^2\,L^2=0. 
\mytag{1.6}
$$
The other variables $E_{20}$, $E_{30}$, $E_{02}$, $E_{03}$, $E_{21}$, 
$E_{12}$ are expressed as rational functions of $E_{10}$, $E_{01}$,
$E_{11}$, and $L$ (see formulas \thetag{4.1}, \thetag{4.3}, \thetag{5.1}, 
\thetag{5.2}, \thetag{4.6}, \thetag{4.7} in \mycite{48}).\par
     The equation \mythetag{1.6} was solved by John Ramsden in \mycite{70}. 
In the case of a rational perfect cuboid, where $L=1$, omitting some inessential
special cases, the general solution of the equation \mythetag{1.6} is given
by the formulas 
$$
\gather
\hskip -2em
E_{11}=-\frac{b\,(c^2+2-4\,c)}{b^2\,c^2+2\,b^2-3\,b^2\,c+c-b\,c^2\,+2\,b},
\mytag{1.7}\\
\vspace{2ex}
\hskip -2em
E_{10}=-\frac{b^2\,c^2+2\,b^2-3\,b^2\,c\,-c}{b^2\,c^2+2\,b^2-3\,b^2\,c
+c-b\,c^2+2\,b},
	\mytag{1.8}\\
\vspace{2ex}
\hskip -2em
E_{01}=-\frac{b\,(c^2+2-2\,c)}{b^2\,c^2+2\,b^2-3\,b^2\,c+c-b\,c^2+2\,b}.
\mytag{1.9}
\endgather
$$
Below are the formulas for $E_{12}$, $E_{21}$, $E_{03}$, $E_{30}$, $E_{02}$, 
$E_{20}$ in \mythetag{1.3}, \mythetag{1.4}, and \mythetag{1.5}:
$$
\allowdisplaybreaks
\gather
\hskip -2em
\gathered
E_{12}=(16\,b^6+32\,b^5-6\,c^5\,b^2+2\,c^5\,b-62\,b^5\,c^6
+62\,b^6\,c^6+16\,b^4\,-\\
-\,180\,b^6\,c^5-c^7\,b^3+18\,b^5\,c^7-12\,b^6\,c^7-2\,b^5\,c^8
+b^6\,c^8+248\,b^5\,c^2\,+\\
+\,248\,b^6\,c^2-96\,b^6\,c+321\,b^6\,c^4-180\,b^5\,c^3-144\,b^5\,c
-360\,b^6\,c^3\,+\\
+\,b^4\,c^8+8\,b^4\,c^6-6\,b^4\,c^7+18\,b^4\,c^5+7\,b^3\,c^6
+90\,b^5\,c^5-14\,b^3\,c^5\,+\\
+\,17\,b^2\,c^4+32\,b^4\,c^2+28\,b^3\,c^3-28\,b^3\,c^2-4\,b\,c^3+8\,b^3\,c
-57\,b^4\,c^4\,+\\
+\,36\,b^4\,c^3-12\,b^2\,c^3-48\,b^4\,c-c^4)\,(b^2\,c^4-6\,b^2\,c^3
+13\,b^2\,c^2\,-\\
-\,12\,b^2\,c+4\,b^2+c^2)^{-1}\,(b\,c-1-b)^{-2}\,(b\,c-c-2\,b)^{-2},
\endgathered\qquad
\mytag{1.10}\\
\vspace{2ex}
\gathered
E_{21}=\frac{b}{2}\,(5\,c^6\,b-2\,c^6\,b^2+52\,c^5\,b^2-16\,c^5\,b
-2\,c^7\,b^2+2\,b^4\,c^8\,-\\
-\,26\,b^4\,c^7-426\,b^4\,c^5-61\,b^3\,c^6+100\,b^3\,c^5
+14\,c^7\,b^3-c^8\,b^3-20\,b\,c^2\,-\\
-\,8\,b^2\,c^2-16\,b^2\,c-128\,b^2\,c^4-200\,b^3\,c^3
+244\,b^3\,c^2+32\,b\,c^3\,+\\
+\,768\,b^4\,c^4-852\,b^4\,c^3+568\,b^4\,c^2+104\,b^2\,c^3-208\,b^4\,c
+8\,c^4\,+\\
+16\,b^3-112\,b^3\,c+142\,b^4\,c^6
+32\,b^4-2\,c^5)\,(b^2\,c^4-6\,b^2\,c^3+13\,b^2\,c^2\,-\\
-12\,b^2\,c-4\,c^3+4\,b^2+c^2)^{-1}\,(b\,c-1-b)^{-2}\,(b\,c-c-2\,b)^{-2},
\endgathered\qquad\quad
\mytag{1.11}\\
\vspace{1ex}
\hskip -2em
\gathered
E_{03}=\frac{b}{2}\,(b^2\,c^4-5\,b^2\,c^3+10\,b^2\,c^2-10\,b^2\,c+4\,b^2
+2\,b\,c+2\,c^2\,-\\
-\,b\,c^3)\,(2\,b^2\,c^4-12\,b^2\,c^3+26\,b^2\,c^2-24\,b^2\,c
+\,8\,b^2-c^4\,b+3\,b\,c^3\,-\\
-\,6\,b\,c+4\,b+c^3-2\,c^2+2\,c)\,(b^2\,c^4-6\,b^2\,c^3+13\,b^2\,c^2\,-\\
-12\,b^2\,c+4\,b^2+c^2)^{-1}\,(b\,c-1-b)^{-2}\,(-c+b\,c-2\,b)^{-2},
\endgathered\qquad\quad
\mytag{1.12}\\
\vspace{1ex}
\hskip -2em
\gathered
E_{30}=c\,b^2\,(1-c)\,(c-2)\,(b\,c^2-4\,b\,c+2+4\,b)
\,(2\,b\,c^2-c^2-4\,b\,c\,+\\
+\,2\,b)\,(b^2\,c^4-6\,b^2\,c^3+13\,b^2\,c^2-12\,b^2\,c+4\,b^2
+c^2)^{-1}\,\times\\
\times\,(b\,c-1-b)^{-2}\,(-c+b\,c-2\,b)^{-2},
\endgathered\qquad\quad
\mytag{1.13}\\
\vspace{1ex}
\hskip -2em
\gathered
E_{02}=\frac{1}{2}\,(28\,b^2\,c^2-16\,b^2\,c-2\,c^2-4\,b^2-b^2\,c^4
+4\,b^3\,c^4-12\,b^3\,c^3\,+\\
+\,4\,b\,c^3+24\,b^3\,c-8\,b\,c-2\,b^4\,c^4+12\,b^4\,c^3-26\,b^4\,c^2
-8\,b^2\,c^3\,+\\
+24\,b^4\,c-16\,b^3-8\,b^4)\,(b\,c-1-b)^{-2}\,(b\,c-c-2\,b)^{-2},
\endgathered\qquad\quad
\mytag{1.14}\\
\vspace{1ex}
\hskip -2em
\gathered
E_{20}=\frac{b}{2}\,(b\,c^2-2\,c-2\,b)\,(2\,b\,c^2-c^2-6\,b\,c+2
+4\,b)\,\times\\
\times\,(b\,c-1-b)^{-2}\,(b\,c-c-2\,b)^{-2}.
\endgathered\qquad\quad
\mytag{1.15}
\endgather
$$
The formulas \mythetag{1.10}, \mythetag{1.11}, \mythetag{1.12}, 
\mythetag{1.13}, \mythetag{1.14}, \mythetag{1.15} were derived in 
\mycite{71} by substituting the formulas \mythetag{1.7}, \mythetag{1.8}, 
and \mythetag{1.9} along with $L=1$ into the corresponding formulas 
from \mycite{48}.\par
     Thus, the right hand sides of the equalities \mythetag{1.3},
\mythetag{1.4}, and \mythetag{1.5} turned out to be expressed through
two arbitrary rational parameters $b$ and $c$. The next step was to
resolve these equalities with respect to $x_1$, $x_2$, $x_3$, $d_1$, 
$d_2$, $d_3$. For this purpose in \mycite{71} the following two cubic 
equations were written:
$$
\align
&\hskip -2em
x^3-E_{10}\,x^2+E_{20}\,x-E_{30}=0,
\mytag{1.16}\\
\vspace{1ex}
&\hskip -2em
d^{\kern 1pt 3}-E_{01}\,d^{\kern 1pt 2}+E_{02}\,d-E_{03}=0.
\mytag{1.17}
\endalign
$$
Note that the left hand sides of the equalities \mythetag{1.3} are regular
symmetric polynomials of the variables $x_1$, $x_2$, $x_3$ (see \mycite{72}). 
Similarly, the left hand sides of the equalities \mythetag{1.4} are regular
symmetric polynomials of the variables $d_1$, $d_2$, $d_3$. For this reason
$x_1$, $x_2$, $x_3$ can be found as roots of the cubic equation 
\mythetag{1.16}. Similarly, $d_1$, $d_2$, $d_3$ are roots of the second
cubic equation \mythetag{1.17}. Relying on these facts, in \mycite{71} the
following two inverse problems were formulated. 
\myproblem{1.1} Find all pairs of rational numbers $b$ and $c$ for which the
cubic equations \mythetag{1.16} and \mythetag{1.17} with the coefficients given
by the formulas \mythetag{1.8}, \mythetag{1.9},	\mythetag{1.12}, \mythetag{1.13}, 
\mythetag{1.14},	\mythetag{1.15} possess positive rational roots $x_1$, $x_2$, 
$x_3$, $d_1$, $d_2$, $d_3$ obeying the auxiliary polynomial equations 
\mythetag{1.5} whose right hand sides are given by the formulas \mythetag{1.7}, 
\mythetag{1.10}, and \mythetag{1.11}. 
\endproclaim
\myproblem{1.2} Find at least one pair of rational numbers $b$ and $c$ for which 
the cubic equations \mythetag{1.16} and \mythetag{1.17} with the coefficients given
by the formulas \mythetag{1.8}, \mythetag{1.9},	\mythetag{1.12}, \mythetag{1.13}, 
\mythetag{1.14},	\mythetag{1.15} possess positive rational roots $x_1$, $x_2$, 
$x_3$, $d_1$, $d_2$, $d_3$ obeying the auxiliary polynomial equations 
\mythetag{1.5} whose right hand sides are given by the formulas \mythetag{1.7}, 
\mythetag{1.10}, and \mythetag{1.11}.
\endproclaim
     Due to Theorem~\mythetheorem{1.1} the inverse problems~\mytheproblem{1.1} 
and \mytheproblem{1.2} are equivalent to finding all rational perfect cuboids and 
to finding at least one rational perfect cuboid respectively. Singularities of 
the inverse problems~\mytheproblem{1.1} and \mytheproblem{1.2} due to the 
denominators in the formulas \mythetag{1.7} through \mythetag{1.15} were studied
in \mycite{73}. Some special cases where the equations \mythetag{1.3},
\mythetag{1.4}, \mythetag{1.5} are solvable with respect to the cuboid variables
$x_1$, $x_2$, $x_3$ and $d_1$, $d_2$, $d_3$ were found in \mycite{74}. 
However, none of these special cases have produced a perfect cuboid since the 
inequalities 
$$
\pagebreak
\xalignat 3
&\hskip -2em
x_1>0, &&x_2>0, &&x_3>0,\\
\vspace{-1.7ex}
\mytag{1.18}\\
\vspace{-1.7ex}
&\hskip -2em
d_1>0, &&d_2>0, &&d_3>0
\endxalignat 
$$
required for solving the problems~\mytheproblem{1.1} and \mytheproblem{1.2} are
not fulfilled in these special cases.\par
     Again, neglecting the inequalities \mythetag{1.18}, an approach to solving
the equations \mythetag{1.3}, \mythetag{1.4}, \mythetag{1.5} was found in 
\mycite{75}. It exploits the following lemma. 
\mylemma{1.1} A reduced cubic equation $y^3+y^2+D=0$ has three rational roots 
if and only if there is a rational number $w$ satisfying the sextic equation 
$$
\hskip -2em
D\,(w^2+3)^3+4\,(w-1)^2\,(1+w)^2=0.
\mytag{1.19}
$$
In this case the roots of the cubic equation $y^3+y^2+D=0$ are given by the
formulas 
$$
\xalignat 3
&y_1=-\frac{2\,(w+1)}{w^2+3},   
&&y_2=\frac{2\,(w-1)}{w^2+3},
&&y_3=\frac{1-w^2}{w^2+3}.
\endxalignat
$$
\endproclaim
     The idea of Lemma~\mythelemma{1.1} belongs to John Ramsden. Its detailed 
proof is given in \mycite{75}, where this lemma was applied to \mythetag{1.16} and 
\mythetag{1.17} (see Lemma 2.1 in \mycite{75}). As a result two sextic equations 
of the form \mythetag{1.19} were derived:
$$
\align
&\hskip -2em
D_1\,(w_1^2+3)^3+4\,(w_1-1)^2\,(1+w_1)^2=0,
\mytag{1.20}\\
\vspace{1ex}
&\hskip -2em
D_2\,(w_2^2+3)^3+4\,(w_2-1)^2\,(1+w_2)^2=0.
\mytag{1.21}\\
\endalign
$$
The $D$-parameters $D_1$ and $D_2$ of the sextic equations \mythetag{1.20} and
\mythetag{1.21} depend on the same two rational numbers $b$ and $c$ as $E_{11}$,  
$E_{10}$, $E_{01}$, $E_{12}$, $E_{21}$, $E_{03}$, $E_{30}$, $E_{02}$, $E_{20}$ 
in the formulas \mythetag{1.7} through \mythetag{1.15}. Therefore the equations
\mythetag{1.20} and \mythetag{1.21} define two algebraic surfaces or, which is
equivalent, two algebraic functions $w_1(b,c)$ and $w_2(b,c)$. Here are the
explicit formulas for the $D$-parameters $D_1$ and $D_2$ in the sextic equations 
\mythetag{1.20} and \mythetag{1.21}:
$$
\gathered
D_1=-\frac{2}{27}\,(7812\,b^4\,c^4\,-216\,b^2\,c^4-52\,b^2\,c^3+1764\,b^3\,c^4
-1200\,b^4\,c^3\,-\\
-\,1848\,b^4\,c^2+720\,b^4\,c-36\,c^4\,b-1512\,b^3\,c^3-36\,c^8\,b^3
+288\,b^3\,c^2\,-\\
-\,108\,c^6\,b^2+380\,c^5\,b^2+378\,c^7\,b^3-231\,c^8\,b^4-300\,c^7\,b^4
+3906\,c^6\,b^4\,-\\
-13\,c^7\,b^2-8904\,c^5\,b^4-882\,c^6\,b^3+18\,c^6\,b-1319\,b^6\,c^8
+20952\,b^5\,c^3\,-\\
-\,11952\,b^5\,c^2+2592\,b^5\,c-48372\,b^6\,c^4+31620\,b^6\,c^3-10552\,b^6\,c^2\,+\\
+\,\,816\,b^6\,c+1494\,b^5\,c^8-5238\,b^5\,c^7-4\,c^5+7905\,b^6\,c^7
-24186\,b^6\,c^6\,+\\
+\,288\,b^6+43740\,b^6\,c^5+7686\,b^5\,c^6+576\,b^7+128\,b^8-15372\,b^5\,c^4\,-\\
-\,1080\,b^7\,c^8-3546\,b^7\,c^6+51\,c^9\,b^6+400\,b^8\,c^8-162\,c^9\,b^5
+8640\,b^7\,c^2\,-\\
-\,3456\,b^7\,c+2808\,b^7\,c^7-1560\,b^8\,c^7+3940\,b^8\,c^6+216\,c^9\,b^7
-960\,b^8\,c\,-\\
-\,6240\,b^8\,c^3+9\,c^{10}\,b^6+7880\,b^8\,c^4+4\,c^{10}\,b^8-6732\,b^8\,c^5 
+45\,c^9\,b^4\,+\\
+\,3200\,b^8\,c^2-11232\,b^7\,c^3+7092\,b^7\,c^4-18\,c^{10}\,b^7
-60\,c^9\,b^8)^2\,(2\,c^2\,+\\
+\,2\,b^4\,c^4-12\,b^4\,c^3+26\,b^4\,c^2-24\,b^4\,c+8\,b^4-6\,b^3\,c^4
+18\,b^3\,c^3\,-\\
-\,36\,b^3\,c+24\,b^3+3\,b^2\,c^4+8\,b^2\,c^3-36\,b^2\,c^2+16\,b^2\,c+12\,b^2
-6\,b\,c^3\,+\\
+\,12\,b\,c)^{-3}\,(b^2\,c^4-6\,b^2\,c^{-3}+13\,b^2\,c^2-12\,b^2\,c+4\,b^2+c^2)^{-2},
\endgathered
\mytag{1.22}
$$
% continued on the next page
$$
\gathered
D_2=-\frac{2\,b^2}{27}\,(832\,b^2\,c^2-1440\,b^2\,c^4-840\,b^2\,c^3
+4788\,b^3\,c^4+396\,b\,c^3\,+\\
+\,720\,b^3\,c+808\,b^4\,c^4+3032\,b^4\,c^3-2576\,b^4\,c^2
-96\,b^4\,c+448\,b^4\,-\\
-\,504\,c^4\,b-4176\,b^3\,c^3-9\,c^8\,b^3+72\,b^3\,c^2
-720\,c^6\,b^2+2288\,c^5\,b^2\,+\\
+\,1044\,c^7\,b^3-322\,c^8\,b^4+758\,c^7\,b^4+404\,c^6\,b^4
-210\,c^7\,b^2-2464\,c^5\,b^4\,-\\
-\,2394\,c^6\,b^3+72\,c^4+252\,c^6\,b+3168\,b^6\,c^8+441\,c^9\,b^5
-7056\,b^5\,c\,+\\
+\,57960\,b^6\,c^4-47232\,b^6\,c^3+25344\,b^6\,c^2
-8064\,b^6\,c-1809\,b^5\,c^8\,+\\
+\,14472\,b^5\,c^2+3951\,b^5\,c^7-72\,c^5+36\,c^6-11808\,b^6\,c^7
+1440\,b^5\,+\\
+\,28980\,b^6\,c^6-49032\,b^6\,c^5-4410\,b^5\,c^6
+8820\,b^5\,c^4-15804\,b^5\,c^3\,+\\
+\,1152\,b^6-504\,c^9\,b^6-45\,c^9\,b^3-6\,c^9\,b^4+104\,c^8\,b^2
+36\,c^{10}\,b^6\,+\\
+\,14\,c^{10}\,b^4-45\,c^{10}\,b^5-99\,c^7\,b)^2\,(6\,b^4\,c^4
-36\,b^4\,c^3+78\,b^4\,c^2
-72\,b^4\,c\,+\\
+\,24\,b^4-12\,b^3\,c^4+36\,b^3\,c^3-72\,b^3\,c+48\,b^3+5\,b^2\,c^4
+16\,b^2\,c^3\,-\\
-\,68\,b^2\,c^2+32\,b^2\,c+20\,b^2-12\,b\,c^3+24\,b\,c
+6\,c^2)^{-3}\,(b^2\,c^4-6\,b^2\,c^3\,+\\
+\,13\,b^2\,c^2-12\,b^2\,c+4\,b^2+c^2)^{-2}.
\endgathered\quad
\mytag{1.23}
$$
The main result of \mycite{75} is expressed by the following two theorems.
\mytheorem{1.2} Each rational point of the algebraic surface \mythetag{1.20}, 
except for points belonging to some algebraic subvariety of this surface, 
determine six rational numbers $x_1$, $x_2$, $x_3$ and $d_1$, $d_2$, $d_3$ 
obeying the cuboid factor equations \mythetag{1.2} as well as the equations 
\mythetag{1.3}, \mythetag{1.4}, \mythetag{1.5} whose right hand sides are given 
by the formulas \mythetag{1.7} through \mythetag{1.15}. 
\endproclaim
\mytheorem{1.3} Each rational point of the algebraic surface \mythetag{1.21}, 
except for points belonging to some algebraic subvariety of this surface, 
determine six rational numbers $x_1$, $x_2$, $x_3$ and $d_1$, $d_2$, $d_3$ 
obeying the cuboid factor equations \mythetag{1.2} as well as the equations 
\mythetag{1.3}, \mythetag{1.4}, \mythetag{1.5} whose right hand sides are given 
by the formulas \mythetag{1.7} through \mythetag{1.15}. 
\endproclaim
     The rational numbers $x_1$, $x_2$, $x_3$ and $d_1$, $d_2$, $d_3$ in 
Theorem~\mythetheorem{1.2} are given by explicit formulas expressing them as 
rational functions of $b$, $c$ and $w_1$:
$$
\xalignat 2
&\hskip -2em
x_i=x_i(b,c,w_1),
&&d_i=d_i(b,c,w_1).
\mytag{1.24}
\endxalignat
$$
But the formulas for $x_i(b,c,w_1)$ and $d_i(b,c,w_1)$ in \mythetag{1.24} are 
very huge. For this reason they are nor presented here. They are given in a 
machine readable form in the ancillary file 
\myhref{http://arxiv.org/src/1209.0723/anc/Solution\podcherkivanie 1.txt}
{{\bf Solution\_\kern 1.5pt 1.txt}} to \mycite{76}. As for the exceptional
subvariety mentioned in Theorem~\mythetheorem{1.2}, it is determined by
the denominator of \mythetag{1.22} and by the denominators in those huge
formulas for $x_1$, $x_2$, $x_3$ and $d_1$, $d_2$, $d_3$.\par
     The rational numbers $x_1$, $x_2$, $x_3$ and $d_1$, $d_2$, $d_3$ in 
Theorem~\mythetheorem{1.3} are also given by explicit formulas expressing them 
as rational functions of $b$, $c$ and $w_2$:
$$
\xalignat 2
&\hskip -2em
x_i=x_i(b,c,w_2),
&&d_i=d_i(b,c,w_2).
\mytag{1.25}
\endxalignat
$$
These formulas for $x_i(b,c,w_2)$ and $d_i(b,c,w_2)$ in \mythetag{1.25} are 
also very huge. For this reason we do not present them here. They are given in 
a machine readable form in the ancillary file \myhref{http://arxiv.org/src/1209.0723/anc/Solution\podcherkivanie 2.txt}
{{\bf Solution\_\kern 1.5pt 2.txt}} to \mycite{76}. The exceptional subvariety 
in Theorem~\mythetheorem{1.3} is determined by the denominator of \mythetag{1.23} 
and by the denominators in those huge formulas for $x_1$, $x_2$, $x_3$ and 
$d_1$, $d_2$, $d_3$.\par
      The algebraic surfaces \mythetag{1.20} and \mythetag{1.21} are not 
independent. They were studied in \mycite{76}. As it was shown in \mycite{76}, 
these two algebraic surfaces are birationally equivalent. This birational equivalence 
is established by two rational functions
$$
\xalignat 2
&\hskip -2em
w_2=w_2(b,c,w_1),
&&w_1=w_1(b,c,w_2).
\mytag{1.26}
\endxalignat
$$
The formulas \mythetag{1.26} should be treated modulo the equations \mythetag{1.20}
and \mythetag{1.21} respectively. Then they produce two mutually inverse 
transformations.\par 
     The functions \mythetag{1.24} and \mythetag{1.25} are also not independent. 
They are related to each other through the birational transformations \mythetag{1.26}:
$$
\align
&\hskip -1em
\aligned
&x_i(b,c,w_1)=x_i(b,c,w_2(b,c,w_1)),\\
&d_i(b,c,w_1)=d_i(b,c,w_2(b,c,w_1)),
\endaligned
\mytag{1.27}\\
\vspace{1ex}
&\hskip -1em
\aligned
&x_i(b,c,w_2)=x_i(b,c,w_1(b,c,w_2)),\\
&d_i(b,c,w_2)=d_i(b,c,w_1(b,c,w_2)).
\endaligned
\mytag{1.28}
\endalign
$$
The explicit formulas for the rational functions $w_2(b,c,w_1)$ and $w_1(b,c,w_2)$
in \mythetag{1.26} are very huge. They are not presented here, but they are given
in a machine readable form in the ancillary file 
\myhref{http://arxiv.org/src/1209.0723/anc/Conversion\podcherkivanie formulas.txt}
{{\bf Conversion\_\kern 1.3pt formulas.txt}} to \mycite{76}.\par
     Due to \mythetag{1.26}, \mythetag{1.27}, and \mythetag{1.28} the algebraic 
surfaces given by the sextic equations \mythetag{1.20} and \mythetag{1.21} are 
closely related to each other. Below we study both of them. Since $D_1=D_1(b,c)$
and $D_2=D_2(b,c)$ in \mythetag{1.20} and \mythetag{1.21}, they are surfaces in
$\Bbb R^3$. The main goal of the present paper is to embed these surfaces into 
$\Bbb R^4$ and show that each point of any one of them lies in the intersection 
of some genus zero curve and some genus one curve specific to this point.\par
     Note that some genus one curves associated with perfect cuboids were 
considered in \mycite{77}. However, they correspond to some very special solutions 
of the cubic equations \mythetag{1.16} and \mythetag{1.17}. Genus one curves 
considered in this paper cover the general case in the theory of the cuboid factor 
equations.\par
\head
2. The structure of $D_1$ and $D_2$ and the 
{\twelvesmmi\bigalpha}\kern 0.5pt-parameters. 
\endhead
     Let's consider the formulas \mythetag{1.22} and \mythetag{1.23} for $D_1$ and $D_2$. 
Looking at these formulas, one can detect the following structure of the expressions for 
$D_1$ and $D_2$:
$$
\xalignat 2
&\hskip -2em
D_1=-\frac{(P_1)^2}{(Q_1)^3},
&&D_2=-\frac{(P_2)^2}{(Q_2)^3}.
\mytag{2.1}
\endxalignat
$$
The denominators $Q_1$ and $Q_2$ in \mythetag{2.1} are given by the formulas
$$
\allowdisplaybreaks
\gather
\hskip -2em
\gathered
Q_1=\frac{3}{2}\,(2\,c^2+2\,b^4\,c^4-12\,b^4\,c^3+26\,b^4\,c^2-24\,b^4\,c+8\,b^4\,-\\
\vspace{0.5ex}
-\,6\,b^3\,c^4+18\,b^3\,c^3-36\,b^3\,c+24\,b^3+3\,b^2\,c^4+8\,b^2\,c^3\,-\\
-36\,b^2\,c^2+16\,b^2\,c+12\,b^2-6\,b\,c^3+12\,b\,c),
\endgathered
\mytag{2.2}\\
\vspace{2ex}
\hskip -2em
\gathered
Q_2=\frac{3}{2}\,(6\,b^4\,c^4-36\,b^4\,c^3+78\,b^4\,c^2-72\,b^4\,c+24\,b^4
-12\,b^3\,c^4\,+\\
\vspace{0.5ex}
+\,36\,b^3\,c^3-72\,b^3\,c+48\,b^3+5\,b^2\,c^4+16\,b^2\,c^3-68\,b^2\,c^2\,+\\
+\,32\,b^2\,c+20\,b^2-12\,b\,c^3+24\,b\,c+6\,c^2).
\endgathered
\mytag{2.3}
\endgather
$$
The numerators $P_{\kern 1pt 1}$ and $P_{\kern 1pt 2}$ in \mythetag{2.1} 
are given by similar formulas:
$$
\gather
\hskip -2em
\gathered
P_{\kern 1pt 1}=\frac{1}{2}\,(7812\,b^4\,c^4\,-216\,b^2\,c^4-52\,b^2\,c^3
+1764\,b^3\,c^4-1200\,b^4\,c^3\,-\\
-\,1848\,b^4\,c^2+720\,b^4\,c-36\,c^4\,b-1512\,b^3\,c^3-36\,c^8\,b^3
+288\,b^3\,c^2\,-\\
-\,108\,c^6\,b^2+380\,c^5\,b^2+378\,c^7\,b^3-231\,c^8\,b^4-300\,c^7\,b^4
+3906\,c^6\,b^4\,-\\
-13\,c^7\,b^2-8904\,c^5\,b^4-882\,c^6\,b^3+18\,c^6\,b-1319\,b^6\,c^8
+20952\,b^5\,c^3\,-\\
-\,11952\,b^5\,c^2+2592\,b^5\,c-48372\,b^6\,c^4+31620\,b^6\,c^3-10552\,b^6\,c^2\,+\\
+\,\,816\,b^6\,c+1494\,b^5\,c^8-5238\,b^5\,c^7-4\,c^5+7905\,b^6\,c^7
-24186\,b^6\,c^6\,+\\
+\,288\,b^6+43740\,b^6\,c^5+7686\,b^5\,c^6+576\,b^7+128\,b^8-15372\,b^5\,c^4\,-\\
-\,1080\,b^7\,c^8-3546\,b^7\,c^6+51\,c^9\,b^6+400\,b^8\,c^8-162\,c^9\,b^5
+8640\,b^7\,c^2\,-\\
-\,3456\,b^7\,c+2808\,b^7\,c^7-1560\,b^8\,c^7+3940\,b^8\,c^6+216\,c^9\,b^7
-960\,b^8\,c\,-\\
-\,6240\,b^8\,c^3+9\,c^{10}\,b^6+7880\,b^8\,c^4+4\,c^{10}\,b^8-6732\,b^8\,c^5 
+45\,c^9\,b^4\,+\\
+\,3200\,b^8\,c^2-11232\,b^7\,c^3+7092\,b^7\,c^4-18\,c^{10}\,b^7
-60\,c^9\,b^8)\,\times\\
\times\,
(b^2\,c^4-6\,b^2\,c^{-3}+13\,b^2\,c^2-12\,b^2\,c+4\,b^2+c^2)^{-1}\!,
\endgathered
\mytag{2.4}\\
\vspace{2ex}
\hskip -2em
\gathered
P_{\kern 1pt 2}=\frac{b}{2}\,(832\,b^2\,c^2-1440\,b^2\,c^4-840\,b^2\,c^3
+4788\,b^3\,c^4+396\,b\,c^3\,+\\
+\,720\,b^3\,c+808\,b^4\,c^4+3032\,b^4\,c^3-2576\,b^4\,c^2
-96\,b^4\,c+448\,b^4\,-\\
-\,504\,c^4\,b-4176\,b^3\,c^3-9\,c^8\,b^3+72\,b^3\,c^2
-720\,c^6\,b^2+2288\,c^5\,b^2\,+\\
+\,1044\,c^7\,b^3-322\,c^8\,b^4+758\,c^7\,b^4+404\,c^6\,b^4
-210\,c^7\,b^2-2464\,c^5\,b^4\,-\\
-\,2394\,c^6\,b^3+72\,c^4+252\,c^6\,b+3168\,b^6\,c^8+441\,c^9\,b^5
-7056\,b^5\,c\,+\\
+\,57960\,b^6\,c^4-47232\,b^6\,c^3+25344\,b^6\,c^2
-8064\,b^6\,c-1809\,b^5\,c^8\,+\\
+\,14472\,b^5\,c^2+3951\,b^5\,c^7-72\,c^5+36\,c^6-11808\,b^6\,c^7
+1440\,b^5\,+\\
+\,28980\,b^6\,c^6-49032\,b^6\,c^5-4410\,b^5\,c^6
+8820\,b^5\,c^4-15804\,b^5\,c^3\,+\\
+\,1152\,b^6-504\,c^9\,b^6-45\,c^9\,b^3-6\,c^9\,b^4+104\,c^8\,b^2
+36\,c^{10}\,b^6\,+\\
+\,14\,c^{10}\,b^4-45\,c^{10}\,b^5-99\,c^7\,b)\,
(b^2\,c^4-6\,b^2\,c^3\,+\\
+\,13\,b^2\,c^2-12\,b^2\,c+4\,b^2+c^2)^{-2}\!.
\endgathered
\mytag{2.5}
\endgather
$$
Due to \mythetag{2.1} and the formulas \mythetag{2.2}, \mythetag{2.3}, 
\mythetag{2.4}, \mythetag{2.5} we can write the sextic equations 
\mythetag{1.20} and \mythetag{1.21} in the following way:
$$
\xalignat 2
&\hskip -2em
\left(\frac{w_1^2+3}{Q_1}\right)^{\lower 3pt\hbox{$\ssize 3$}}
\!=\left(\frac{2\,(w_1^2-1)}{P_{\kern 1pt 1}}\right)^{\lower 3pt
\hbox{$\ssize 2$}}\!\!,
&&\left(\frac{w_2^2+3}{Q_2}\right)^{\lower 3pt\hbox{$\ssize 3$}}
\!=\left(\frac{2\,(w_2^2-1)}{P_{\kern 1pt 2}}\right)^{\lower 3pt
\hbox{$\ssize 2$}}\!\!.\qquad
\mytag{2.6}\\
\endxalignat
$$
\mylemma{2.1} If\/ $x$ and $y$ are two rational numbers obeying the 
equality $x^3=y^2$, then there is a third rational number $\alpha$ such 
that \pagebreak $x=\alpha^2$ and $y=\alpha^3$.
\endproclaim
\demo{Proof} Each rational number is presented as an irreducible ratio of two
integer numbers. Expanding these numbers into prime factors we find that each 
rational number is a product of distinct prime numbers to some integer powers
which can be either positive or negative. For $x$ and $y$ this fact yields
$$
\xalignat 2
&\hskip -2em
x=\pm\,p_1^{\kern 1pt\beta_1}\!\cdot\ldots\cdot p_m^{\kern 1pt\beta_m},
&&y=\pm\,q_1^{\gamma_1}\!\cdot\ldots\cdot q_n^{\gamma_n}.
\mytag{2.7}
\endxalignat
$$
Substituting \mythetag{2.7} into the equality $x^3=y^2$, we find that $x>0$, 
$m=n$ and the prime factors $p_1,\,\ldots,\,p_m$ should coincide with the prime 
factors $q_1,\,\ldots,\,q_n$ up to some permutation. Performing this permutation 
upon the prime factors $q_1,\,\ldots,\,q_n$, for the exponents $\beta_1,\,\ldots,
\,\beta_m$ and $\gamma_1,\,\ldots,\,\gamma_n$ in \mythetag{2.7} we derive 
$$
\hskip -2em
3\,\beta_i=2\,\gamma_i\text{, \ where \ }i=1,\,\ldots,\,n.
\mytag{2.8}
$$
Due to \mythetag{2.8} there are some unique integer numbers $\alpha_1,\,\ldots,\,
\alpha_n$ such that 
$$
\xalignat 2
&\hskip -2em
\beta_i=2\,\alpha_i,
&&\gamma_i=3\,\alpha_i.
\mytag{2.9}
\endxalignat
$$
Using the numbers $\alpha_1,\,\ldots,\,\alpha_n$, we define the number $\alpha$ by
setting
$$
\hskip -2em
\alpha=\pm\,p_1^{\kern 1pt\alpha_1}\!\cdot\ldots\cdot p_n^{\kern 1pt\alpha_n}.
\mytag{2.10}
$$
The sign of $\alpha$ in \mythetag{2.10} is chosen to be coinciding with the sign of 
$y$ in \mythetag{2.7}. Then from \mythetag{2.9}, using $x>0$ and \mythetag{2.7}, we
derive the required equalities $x=\alpha^2$ and $y=\alpha^3$. Lemma~\mythelemma{2.1}
is proved.
\qed\enddemo
     Now we can apply Lemma~\mythelemma{2.1} to the equations \mythetag{2.6}. As a 
result we obtain the following two theorems whose proofs are obvious. 
\mytheorem{2.1} For each rational point $(b,c,w_1)$ of the algebraic surface
\mythetag{1.20} in $\Bbb R^3$, where $w_1\neq\pm\,1$, there is a rational number 
$\alpha_1$ such that 
$$
\xalignat 2
&\hskip -2em
w_1^2+3=Q_1\,\alpha_1^2,
&&2\,(w_1^2-1)=P_{\kern 1pt 1}\,\alpha_1^3.
\mytag{2.11}
\endxalignat
$$
The numbers $P_{\kern 1pt 1}$ and $Q_1$ in \mythetag{2.11} are given by the formulas 
\mythetag{2.2} and \mythetag{2.4}.
\endproclaim
\mytheorem{2.2} For each rational point $(b,c,w_2)$ of the algebraic surface
\mythetag{1.21} in $\Bbb R^3$, where $w_2\neq\pm\,1$, there is a rational number 
$\alpha_2$ such that 
$$
\xalignat 2
&\hskip -2em
w_2^2+3=Q_2\,\alpha_2^2,
&&2\,(w_2^2-1)=P_{\kern 1pt 2}\,\alpha_2^3.
\mytag{2.12}
\endxalignat
$$
The numbers $P_{\kern 1pt 2}$ and $Q_2$ in \mythetag{2.12} are given by the formulas 
\mythetag{2.3} and \mythetag{2.5}.
\endproclaim
\head
3. Genus zero curves.
\endhead
     Let's consider the equations \mythetag{2.11} separately. Since $Q_1=Q_1(b,c)$,
the first equation \mythetag{2.11} defines a hypersurface in $\Bbb R^4$. Similarly,
$P_{\kern 1pt 1}=P_{\kern 1pt 1}(b,c)$. Therefore the second equation \mythetag{2.11}
defines another hypersurface in $\Bbb R^4$. Thus, we get an embedding of the algebraic
surface \mythetag{1.20} into $\Bbb R^4$ where this surface is presented as the 
intersection of two algebraic hypersurfaces. In a similar way, the equations 
\mythetag{2.12} define an embedding of the algebraic surface \mythetag{1.21} into 
$\Bbb R^4$ \pagebreak where this surface is presented as the intersection of two 
algebraic hypersurfaces.\par
     Now assume that the rational numbers $b$ and $c$ are chosen to be constants. Then 
$P_{\kern 1pt 1}$ and $Q_1$ in \mythetag{2.11} are also two rational constants. Under
this assumption the first equation \mythetag{2.11} defines a genus zero curve in 
$\Bbb R^2$, while the second equation defines a genus one curve. The same is true for 
the equations \mythetag{2.12}, provided $b$ and $c$ are considered as constants. 
\mydefinition{3.1} A genus zero algebraic curve is called a rational curve if it is 
birationally equivalent to a line (see \mycite{78}). 
\enddefinition
\mytheorem{3.1} An algebraic curve in $\Bbb R^2$ given by the quadratic equation
$$
\hskip -2em
w^2+3=Q\,\alpha^2,
\mytag{3.1}
$$
where $Q$ is a rational number, is birationally equivalent to a line over $\Bbb Q$
if and only if it has at least one rational point $(w_0,\alpha_0)$ in $\Bbb R^2$. 
\endproclaim
\demo{Proof} Necessity. Assume that the curve \mythetag{3.1} is birationally 
equivalent to a line over $\Bbb Q$. Then there are two rational functions 
with rational coefficients 
$$
\xalignat 2
&\hskip -2em
w=w(t),
&&\alpha=\alpha(t)
\mytag{3.2}
\endxalignat
$$
whose domain is some Zariski open subset $D\subseteq\Bbb Q\subset\Bbb R$ and 
such that they satisfy the equality \mythetag{3.1} identically. Choosing some 
point $t_0\in D$ and substituting $t=t_0$ into \mythetag{3.2}, we get two rational 
numbers $w_0=w(t_0)$ and $\alpha_0=\alpha(t_0)$ satisfying the equality 
\mythetag{3.1}. They constitute a rational point $(w_0,\alpha_0)$ of our curve 
\mythetag{3.1}.\par
     Sufficiency. Assume that $(w_0,\alpha_0)$ is some rational point of the
curve \mythetag{3.1}. Then the following equality holds for $w_0$ and $\alpha_0$:
$$
\hskip -2em
w_0^2+3=Q\,\alpha_0^2.
\mytag{3.3}
$$
Subtracting \mythetag{3.3} from \mythetag{3.1}, we derive
$$
\hskip -2em
(w-w_0)\,(w+w_0)=Q\,(\alpha-\alpha_0)\,(\alpha+\alpha_0).
\mytag{3.4}
$$
The equality \mythetag{3.4} can be transformed to the following one:
$$
\hskip -2em
\frac{w-w_0}{\alpha+\alpha_0}=Q\,\frac{\alpha-\alpha_0}{w+w_0}.
\mytag{3.5}
$$
Using the equality \mythetag{3.5}, we introduce the parameter $t$ by 
setting
$$
\xalignat 2
&\hskip -2em
\frac{\alpha-\alpha_0}{w+w_0}=t,
&&\frac{w-w_0}{\alpha+\alpha_0}=Q\,t.
\mytag{3.6}
\endxalignat
$$
The equations \mythetag{3.6} are easily resolved with respect to $w$
and $\alpha$:
$$
\xalignat 2
&\hskip -2em
w=\frac{w_0+2\,Q\,\alpha_0\,t+Q\,w_0\,t^2}{1-Q\,t^2},
&&\alpha=\frac{\alpha_0+2\,w_0\,t+Q\,\alpha_0\,t^2}{1-Q\,t^2}.
\qquad
\mytag{3.7}
\endxalignat
$$
The functions \mythetag{3.7} are that very rational functions $w(t)$ 
and $\alpha(t)$ which implement the birational equivalence of the curve 
\mythetag{3.1} and a line. Theorem is proved.
\qed\enddemo
     Note that $Q$ in \mythetag{3.1} is a rational number. Therefore
$Q=(M\,m^2)/(N\,n^2)$, where $M$ and $N$ are relatively prime square 
free integer numbers. Let $(w_0,\alpha_0)$ be a rational point of the 
algebraic curve \mythetag{3.1}. Then $w_0$ and $(\alpha_0\,m)/(N\,n)\,$ 
are two rational numbers. We can express these numbers as 
$$
\xalignat 2
&\hskip -2em
w_0=\frac{X}{Z},
&&\frac{\alpha_0\,m}{N\,n}=\frac{Y}{Z},
\mytag{3.8}
\endxalignat
$$
where $X$, $Y$, and $Z$ are three integers. Applying $Q=(M\,m^2)/(N\,n^2)$ 
and \mythetag{3.8} to \mythetag{3.1}, we derive the following Diophantine
equation:
$$
\hskip -2em
X^2-M\,N\,Y^2+3\,Z^2=0.
\mytag{3.9}
$$
The Diophantine equation \mythetag{3.9} is an instance of the Legendre equation
(see \mycite{79}, \mycite{80}, and \mycite{81}). There is a criterion for the
equation \mythetag{3.9} to have a solution. This criterion can be derived from 
the following Legendre theorem.
\mytheoremwithtitle{3.2}{ {\rm(}Legendre{\rm)}} A quadratic Diophantine equation
with square free and pairwise relatively prime coefficients 
$$
\hskip -2em
A\,X^2+B\,Y^2+C\,Z^2=0
\mytag{3.10}
$$
has a nonzero solution if and only if its coefficients are not all of the same 
sign and if\/ $-B\,C$, $-C\,A$, and $-A\,B$ are squares modulo $A$, $B$, and $C$ 
respectively. 
\endproclaim
     Due to \mythetag{3.8} the numbers $M$ and $N$ in \mythetag{3.9} are 
relatively prime and square free. Therefore their product $M\,N$ is also
square free. The equation \mythetag{3.9} is a special instance of the equation
\mythetag{3.10}. There are special solvability criteria for it (see \mycite{80}).
\mytheorem{3.3} If\/ $M\,N$ is nonzero modulo $3$, then the Diophantine
equation \mythetag{3.9} has a nonzero solution if and only if the following
three conditions are fulfilled:
$$
\xalignat 3
&\text{1) }M\,N>0,\quad 
&&\text{2) }-3\text{\ is a square modulo }M\,N,\quad
&&\text{3) }M\,N\text{\ is a square modulo }3.
\endxalignat
$$
\endproclaim
\mytheorem{3.4} If\/ $M\,N$ is zero modulo $3$, then the Diophantine
equation \mythetag{3.9} has a nonzero solution if and only if the following
three conditions are fulfilled:
$$
\xalignat 3
&\text{1) }M\,N>0,\quad 
&&\text{2) }-3\text{\ is a square modulo }M\,N,\quad
&&\text{3) }\frac{M\,N}{3}\text{\ is a square modulo }3.
\endxalignat
$$
\endproclaim
     Applying Theorems~\mythetheorem{3.3} and \mythetheorem{3.4}, we can study 
the rationality of genus zero curves from \mythetag{2.11} and \mythetag{2.12}
for particular numeric values of $b$ and $c$:
$$
\xalignat 2
&\hskip -2em
w_1^2+3=Q_1\,\alpha_1^2,
&&w_2^2+3=Q_2\,\alpha_2^2.
\mytag{3.11}
\endxalignat
$$
The results of such numeric studies for the curves \mythetag{3.11} are presented 
in the ancillary file {\bf Quadratic\_\kern 1.5pt curves.txt} attached to this 
arXiv submission.\par
     Looking through the numeric output file {\bf Quadratic\_\kern 1.5pt 
curves.txt}, we see that some genus zero curves are rational over $\Bbb Q$. 
Others are not rational. So, both types of curves arise in the theory of 
cuboid factor equations.\par
\head
4. Genus one curves.
\endhead
    Genus one algebraic curves in \mythetag{2.11} and \mythetag{2.12} are more 
regular than genus zero curves \mythetag{3.11}. They are given by the following 
cubic equations:
$$
\xalignat 2
&\hskip -2em
2\,(w_1^2-1)=P_{\kern 1pt 1}\,\alpha_1^3,
&&2\,(w_2^2-1)=P_{\kern 1pt 2}\,\alpha_2^3.
\mytag{4.1}
\endxalignat
$$
\mydefinition{4.1} A genus one algebraic curve is called an elliptic curve if it 
has at least one rational point either finite or at infinity (see \mycite{82}). 
\enddefinition
     The parameters $P_{\kern 1pt 1}$ and $P_{\kern 1pt 2}$ in \mythetag{4.1} are
given by the formulas \mythetag{2.4} and \mythetag{2.5}. But despite the values
of these parameters, both curves \mythetag{4.1} have rational points
$$
\xalignat 2
&\hskip -2em
(w_1,\alpha_1)=(\pm\,1,0),
&&(w_2,\alpha_2)=(\pm\,1,0).
\mytag{4.2}
\endxalignat
$$
The points \mythetag{4.2} are exceptional in the sense of Theorems~\mythetheorem{2.1}
and \mythetheorem{2.2}. Nevertheless, they are sufficient to say that both curves 
\mythetag{4.1} are elliptic curves.\par
\head
5. Concluding remarks. 
\endhead
     Summarizing the results of Section 3 and Section 4 above, we conclude that
each rational point of the algebraic surfaces \mythetag{1.20} and \mythetag{1.21}
is produced through the intersection of some elliptic curve and some rational
curve. This fact could be used for to apply powerful tools of the modern theory 
of elliptic curves in studying the equations \mythetag{1.20} and \mythetag{1.21}
and thus in studying perfect cuboids.  

\enddocument
\end